# Waves over a periodic progressive modulation: A `python` tutorial*

Hussein Nassar, Andrew N. Norris and Guoliang Huang

**Abstract** This chapter presents a selection of theoretical and numerical tools suitable for the study of wave propagation in time-dependent media. The focus is on one-dimensional spring-mass chains whose properties are modulated in space and time in a periodic progressive fashion. The chapter is written for the uninitiated newcomer as well as for the theoretically inclined numerical empiricist. Thus, whenever possible, deployed theory is motivated and exploited numerically, and code for example simulations is written in `python`. The chapter begins with an introduction to Mathieu's equation and its stability analysis using the monodromy matrix; generalizations to systems with multiple degrees of freedom are then pursued. The progressive character of the modulation leads to a factorization of the monodromy matrix and provides a "discrete change of variables" otherwise only available for continuous systems. Moreover, the factorization allows to reduce the computational complexity of dispersion diagrams and of long term behaviors. Chosen simulations illustrate salient features of non-reciprocity such as strong left-right biases in the speed and power of propagated waves.

Hussein Nassar
Department of Mechanical and Aerospace Engineering, University of Missouri, Columbia, MO 65211, USA e-mail: `nassarh@missouri.edu`

Andrew N. Norris
Mechanical and Aerospace Engineering, Rutgers University, Piscataway, NJ 08854-8058, USA e-mail: `norris@rutgers.edu`

Guoliang Huang
Department of Mechanical and Aerospace Engineering, University of Missouri, Columbia, MO 65211, USA e-mail: `huangg@missouri.edu`

* Dedicated to the memory of Miles V. Barnhart (1989-2020).

## 1 Introduction

The equation

$$k(u_{n+1} - 2u_n + u_{n-1}) = m\ddot{u}_n \tag{1}$$

governs the propagation of waves along a chain of springs of constant $k$ connecting a series of nodes of mass $m$. Here, $u_n = u_n(t)$ is the displacement of node number $n$ at time $t$. This finite difference equation has been around for a long time; Brillouin reports that Newton first used it in 1686 as a model for sound propagation in air [1]. In that model, $k$ and $m$ are constants. By contrast, the focus here is on cases where $k$ and $m$ are periodically modulated in space and time and, to be more specific, read

$$k = k_n(t) = k(\xi n - \nu t), \quad m = m_n(t) = m(\xi n - \nu t). \tag{2}$$

It is as if properties $k$ and $m$ are now externally driven in a progressive periodic wave-like fashion at a wavenumber $\xi$ and frequency $\nu$ (Fig. 1). The way in which these properties are controlled in practice will not be of concern although it would help the mindful reader to know that such modulations, to certain extent, are experimentally feasible [2]. In any case, the modulation creates a bias: waves traveling with the modulation, that is with a velocity of the same sign as $\nu/\xi$, are expected to behave differently than waves traveling against the modulation, that is with a velocity of the opposite sign as $\nu/\xi$. This bias could be in propagation speeds, transmission and reflection coefficients, and so on. Such media, with a left-right bias, are qualitatively referred to as "non-reciprocal" and can be of use in technological applications.

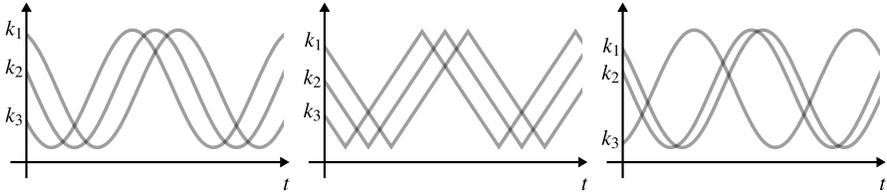

**Fig. 1** Examples of periodic modulations: the left and middle instances are progressive in the sense intended here; the right one is not.

In the modulated chain, the governing equations are

$$\begin{aligned}\dot{u}_n &= p_n/m_n, \\ \dot{p}_n &= k_n(u_{n+1} - u_n) - k_{n-1}(u_n - u_{n-1}).\end{aligned} \tag{3}$$

Thus, the original second-order Ordinary Differential Equation (ODE) is broken into two first-order ODEs by introducing one auxiliary variable: the momentum $p$. One could have chosen velocity $v$ as the auxiliary variable, but the ODE for $\dot{v}$ would have been more involved. The purpose of this chapter is to analyze and exemplify the various behaviors that equation (3) encapsulates using theoretical and numerical tools. Emphasis is on stability issues, or how to determine whether oscillations will grow indefinitely or remain bounded; how to determine the natural modes, their frequencies, and the wavelengths at which they propagate; and on what non-reciprocity looks like and how severe it can be. The main tool of analysis will be the monodromy matrix which is the matrix that propagates initial conditions through a period of time $T = 2\pi/\nu$. The monodromy matrix is first introduced in the context of Mathieu's equation which governs the motion of a single mass; systems with multiple masses arranged in periodic unit cells are subsequently investigated. The main novelty of the present work is a factorization of the monodromy matrix that leverages the progressive character of the modulation. The factorization allows to reduce the computational complexity of dispersion diagrams and of long term behaviors and further provides what is essentially a discrete change of variables otherwise only available for continuous systems.

The chapter is written as a tutorial for the newcomer; it is hopefully self-contained and strikes an appealing balance between analytical rigor and numerical empiricism. This is in particular why an effort has been made to provide code for example simulations written in `python`'s standard libraries `numpy` and `scipy` and why issues of computational complexity and numerical accuracy are occasionally brought up.

## 2 Mathieu's equation

Consider, for starters, the motion of a driven harmonic oscillator without damping:

$$\begin{aligned} \dot{u} &= p/m, \\ \dot{p} &= -ku. \end{aligned} \tag{4}$$

Let, to simplify further, $m$ be a constant and $k = \bar{k} + k'\cos(\nu t)$. It is possible to change variables and substitute $(t/\nu, u/(m\nu))$ for $(t, u)$. In the new variables, the ODEs become

$$\begin{aligned} \dot{u} &= p, \\ \dot{p} &= -(\delta + \epsilon \cos t)u, \end{aligned} \tag{5}$$

with $\delta = \bar{k}/(m\nu^2)$ and $\epsilon = k'/(m\nu^2)$. Thus, parameter $\delta$ measures how slow the modulation is compared to the natural frequency in the absence of modulation, i.e. $\sqrt{\bar{k}/m}$, and parameter $\epsilon$ measures how strong the modulation is. As for time $t$, it is now dimensionless.

Equation (5) is known as Mathieu's equation [3]. In its original context, Mathieu's equation describes the angular components of the free vibration modes of an elliptical membrane. As it turns out, it is also characteristic of parametrically driven oscillators, i.e., oscillators driven by a change in their parameters rather than by a change in the applied force. The main concern here is whether, given a set of initial conditions, oscillations will be periodic, or quasi-periodic, or will grow or decay, and if so, at what rates. Such information can be extracted from the monodromy matrix.

## 2.1 The monodromy matrix

Standard theory of linear ODEs ensures that system (5) has a two-dimensional linear space of solutions. For instance, let

$$\mathbf{\Phi}^i(t) = \begin{bmatrix} u^i(t) \\ p^i(t) \end{bmatrix}, \quad i = 1, 2, \tag{6}$$

be the solutions for the initial conditions $\{u(0) = 1, p(0) = 0\}$ and $\{u(0) = 0, p(0) = 1\}$, respectively. Then, any solution

$$\varphi(t) = \begin{bmatrix} u(t) \\ p(t) \end{bmatrix} \tag{7}$$

is a linear combination

$$\varphi(t) = u(0)\mathbf{\Phi}^1(t) + p(0)\mathbf{\Phi}^2(t) = \mathbf{\Phi}(t)\varphi(0) \tag{8}$$

where $\mathbf{\Phi} = \begin{bmatrix} \mathbf{\Phi}^1 | \mathbf{\Phi}^2 \end{bmatrix}$ is a $2 \times 2$ matrix known as the fundamental matrix of the ODE system. In particular, at time $t = 2\pi$,

$$\varphi(2\pi) = \mathbf{\Phi}(2\pi)\varphi(0) = \mathbf{M}\varphi(0), \tag{9}$$

where the monodromy matrix $\mathbf{M} = \mathbf{\Phi}(2\pi)$ has been introduced. Matrix $\mathbf{M}$ thus propagates a set of initial conditions for a period of time $T = 2\pi$.

Now the coefficients in system (5) are periodic but this does not mean that solutions $\varphi(t)$ are necessarily periodic and, in general, $\varphi(t + 2\pi) \neq \varphi(t)$. What is true nonetheless is that $\varphi(t + 2\pi)$ is also a solution of (5) but with initial conditions $\{u(2\pi), p(2\pi)\}$. Therefore,

$$\varphi(4\pi) = \mathbf{M}\varphi(2\pi) = \mathbf{M}^2\varphi(0). \tag{10}$$

More generally,

$$\varphi(2\pi n) = \mathbf{M}^n\varphi(0). \tag{11}$$

This relationship summarizes the usefulness of the monodromy matrix: to infer the long term behavior of solutions across multiple time periods, it is enough to integrate system (5) over a single time period so as to get $\mathbf{M}$ and then iterate its powers.

Some properties of $\mathbf{M}$ are worth highlighting. On one hand, $\mathbf{M}$ is a real matrix. Thus, if $\lambda$ is a complex eigenvalue, then so is its complex conjugate $\lambda^*$. On the other hand, $\mathbf{M}$ has a determinant of 1. Indeed, letting

$$W(t) = \det \mathbf{\Phi}(t) = u^1(t)p^2(t) - p^1(t)u^2(t), \tag{12}$$

be the Wronskian, it is elementary to check that $\dot{W} = 0$, meaning that $W(t) = W(0) = 1$ and in particular, $\det \mathbf{M} = W(2\pi) = 1$.

## Code: Mathieu's monodromy matrix

The following code is a `python` implementation of a function which computes the monodromy matrix of Mathieu's equation. The code uses `scipy`'s ODE solver `odeint` which employs the "lsoda" method; "lsoda" is of variable order and automatically switches between an Adams scheme, for non-stiff ODEs, and a BDF scheme, for stiff ODEs.[2] In the present context, it seems to prefer an Adams scheme of order 6, 7 or 8.

```python
import numpy as np
from scipy.integrate import odeint

def mathieu(phi, t, delta, epsilon):
    # returns [\dot{u},\dot{p}] for Mathieu's equation
    # with parameters delta, epsilon
    return [phi[1], -(delta + epsilon * np.cos(t)) * phi[0]]

def monodromy(delta, epsilon):
    # returns the monodromy matrix Mon for Mathieu's equation
    # with parameters delta, epsilon and canonical initial conditions
    T = 2 * np.pi
    Mon = np.zeros((2, 2))

    i = 0
    for phi0 in np.eye(2):
        Mon[:, i] = odeint(mathieu, phi0, [0, T],
                           args=(delta, epsilon))[-1, :].T
        i += 1

    return Mon
```

## 2.2 Stability

The oscillator driven with parameters $(\delta, \epsilon)$ is stable if for all initial conditions $\varphi(0)$, the solution $\varphi(t)$ is bounded; if, however, there is a $\varphi(0)$ for which $\varphi(t)$ grows indefinitely, then the oscillator is unstable. It is not hard[3] to see that it is equivalent to require $\varphi(2\pi n)$ be bounded in order to achieve stability. In other words, stability is encoded into the monodromy matrix $\mathbf{M}$.

---

[2] see https://docs.scipy.org/doc/scipy/reference/generated/scipy.integrate.LSODA.html

[3] $\varphi(t)$ bounded $\implies \varphi(2\pi n)$ bounded is immediate. For the reciprocal, one can write $\varphi(t) = \mathbf{\Phi}(t - 2\pi n)\varphi(2\pi n)$ where $0 \le t - 2\pi n \le 2\pi$. Boundedness of $\varphi(t)$ follows from that of $\varphi(2\pi n)$ and from the continuity of $\mathbf{\Phi}$ over $[0, 2\pi]$.

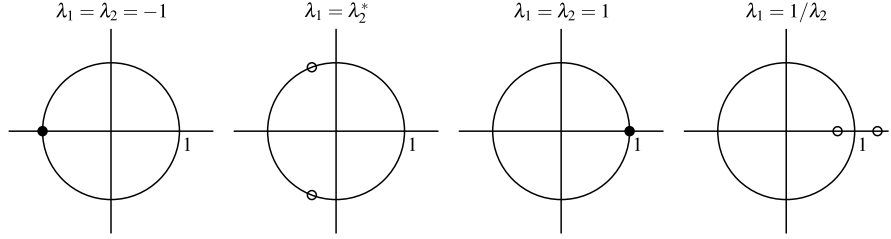

**Fig. 2** Possible placements of the eigenvalues of the monodromy matrix, relative to the unit circle: the first three configurations are stable; the last one is unstable.

Assume $\mathbf{M}$ has two linearly independent eigenvectors, $V^1$ and $V^2$, with eigenvalues $\lambda^1$ and $\lambda^2$. Then $\mathbf{M}$ can be diagonalized into

$$\mathbf{M} = \mathbf{V} \begin{bmatrix} \lambda_1 & 0 \\ 0 & \lambda_2 \end{bmatrix} \mathbf{V}^{-1}, \quad \mathbf{V} = \begin{bmatrix} V^1 | V^2 \end{bmatrix}. \tag{13}$$

Therefore,

$$\mathbf{M}^n = \mathbf{V} \begin{bmatrix} \lambda_1^n & 0 \\ 0 & \lambda_2^n \end{bmatrix} \mathbf{V}^{-1} \tag{14}$$

remains bounded whenever $\lambda^1$ and $\lambda^2$ are on or inside the unit circle. Given that $\mathbf{M}$ is real with a unit determinant, the eigenvalues can be in four qualitatively different configurations all of which are stable except for one; all are depicted on Fig. 2.

Now assume $\mathbf{M}$ has only one eigenvector $V_1$ with a double eigenvalue $\lambda$. Although $\det \mathbf{M} = 1$ implies $\lambda = \pm 1$, the oscillator turns out to be unstable in this case. Indeed, let $V_2$ be any vector independent of $V_1$, then

$$\mathbf{M} = \mathbf{V} \begin{bmatrix} \lambda & \alpha \\ 0 & \lambda \end{bmatrix} \mathbf{V}^{-1}, \quad \alpha \neq 0, \tag{15}$$

and $\mathbf{M}^n V_2 = n\alpha\lambda^{n-1} V_1 + \lambda^n V_2$ grows linearly without bounds.

In conclusion, for the oscillator to transition back and forth between stability and instability, both eigenvalues should coincide, either at 1 or −1. Thus, in the space of the parameters $(\delta, \epsilon)$, transitions occur along curves of equation $\operatorname{tr} \mathbf{M} = \pm 2$. The oscillator is stable when $|\operatorname{tr} \mathbf{M}| < 2$ and is unstable when $|\operatorname{tr} \mathbf{M}| > 2$. When $\operatorname{tr} \mathbf{M} = \pm 2$, the oscillator is stable if $\mathbf{M}$ admits two linearly independent eigenvectors and is unstable otherwise.

### Code: The transition curves of Mathieu's equation

Here is a sample code which draws the transition curves over a range of parameters $(\delta, \epsilon)$. The result is shown on Fig. 3. Note that parameter `Nd` is chosen an order of magnitude larger than `Ne` in order to protect some of the delicate features of the contours tr $\mathbf{M} = \pm 2$ along the $\delta$-axis. Furthermore, it turns out that these contours are rather sensitive to numerical integration errors; thus, the code that produced Fig. 3 employs levels `1.9999` and `-1.9999` instead. This is also why in the function `monodromy`, the "lsoda" method was preferred to the more typical explicit Runge-Kutta of order 5(4) as it resulted in more robust contours.

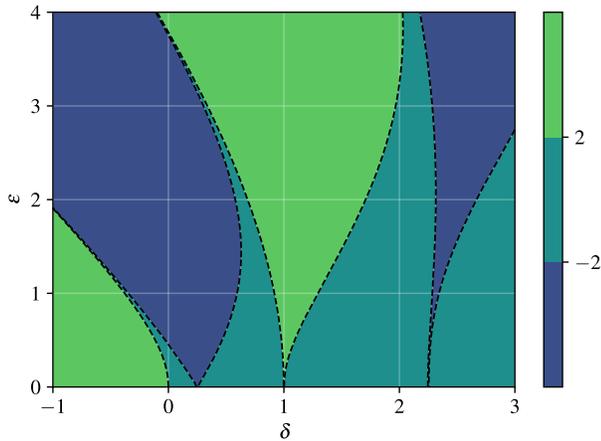

**Fig. 3** Contour plots of tr $\mathbf{M} = \pm 2$ over the $(\delta, \epsilon)$-space: regions where $-2 < \text{tr}\,\mathbf{M} < 2$ are stable; regions where $|\text{tr}\,\mathbf{M}| > 2$ are unstable; crossing the transition curves (dashed lines), the oscillator goes from being stable to unstable, or conversely. The transition curves themselves are unstable, with possible linear growth, except for $\epsilon = 0$. It is customary to highlight that there exists regions of stability even for $\delta < 0$.

```python
# number of samples for delta and epsilon
Nd = 1000
Ne = 100

# initialize plot range and mesh
d_min, d_max = 0, 4
e_min, e_max = 0, 4

[delta, epsilon] = np.meshgrid(np.linspace(d_min, d_max, Nd),
                               np.linspace(e_min, e_max, Ne))

# compute trace over mesh
trMon = np.zeros((Ne, Nd))

for i in range(Nd):
    for j in range(Ne):
        trMon[j, i] = np.trace(monodromy(delta[j, i], epsilon[j, i]))

# plot contours (delta, epsilon, trMon) using favorite library
```

## 3 Periodic progressive modulations

Bringing back space dependence, consider now a chain of $N$ nodes and $N$ springs. The governing first-order ODEs constitute a linear system of $2N$ equations:

$$\dot{u} = M^{-1}p,$$
$$\dot{p} = -Ku,$$
(16)

where $M$ is the mass matrix, $K = A^\top CA$ is the stiffness matrix, $C$ is the constitutive matrix, $A$ is the compatibility matrix and $A^\top$ is its transpose. They read

$$M = \begin{bmatrix} m_1 & & \\ & \ddots & \\ & & m_N \end{bmatrix}, \quad C = \begin{bmatrix} k_1 & & \\ & \ddots & \\ & & k_N \end{bmatrix}, \quad A = \begin{bmatrix} -1 & 1 & & \\ & \ddots & \ddots & \\ & & & 1 \\ 1 & & & -1 \end{bmatrix}, \quad (17)$$

where uncharted territories are filled with 0's. The boundary conditions are taken to be periodic which is why $A$ has a lonely 1 in its bottom left corner.

As for space-time dependence, let

$$k_n = k_o + \delta k \cos(\xi n - \nu t), \quad m_n = m_o + \delta m \cos(\xi n - \nu t), \quad (18)$$

where $k_o$ and $m_o$ are average stiffness and mass values, and $\delta k$ and $\delta m$ are the amplitudes of the perturbation brought by the modulation. Letting $Z$ be the number of nodes per unit cell, the modulation wavenumber reads $\xi = 2\pi/Z$. It is further assumed that the chain contains an integer number of periodic unit cells: $N = Z \times \mathsf{c}$, with $\mathsf{c}$ being the number of periodic unit cells in the chain.

Next, the monodromy matrix is introduced once more. Its eigenvalues and eigenvectors are tapped to build dispersion diagrams which depict relationships between the wavenumbers and frequencies that can propagate through the chain. A factorization of the monodromy matrix leveraging the progressive character of the modulation will permit to shorten numerical integration times as well as to improve the resolution of the dispersion diagrams. Last, a few symptoms of non-reciprocity are explored based on dispersion diagrams and on time-domain simulations. At a certain point, implicit symplectic integrators are briefly introduced to remedy potential numerical instabilities encountered in simulations.

### 3.1 The monodromy matrix

The system of $2N$ ODEs can be recast into the matrix form

$$\dot{\varphi} = \mathbf{A}\varphi, \quad \varphi = \begin{bmatrix} u \\ p \end{bmatrix}, \quad \mathbf{A} = \begin{bmatrix} 0_{N \times N} & M^{-1} \\ -K & 0_{N \times N} \end{bmatrix}. \quad (19)$$

The $2N \times 2N$ fundamental matrix is then solution to the initial value problem

$$\dot{\boldsymbol{\Phi}} = \mathbf{A}\boldsymbol{\Phi}, \quad \boldsymbol{\Phi}(0) = \mathbf{I}, \tag{20}$$

where $\mathbf{I}$ is the identity matrix of order $2N$. Thus, the monodromy matrix is $\mathbf{M} = \boldsymbol{\Phi}(T)$ with $T = 2\pi/\nu$ being the modulation period.

Here too, the monodromy matrix is real and of unit determinant. Indeed, letting $W(t) = \det \boldsymbol{\Phi}(t)$ be the Wronskian, Jacobi's formula yields

$$\dot{W} = \det \boldsymbol{\Phi} \operatorname{tr}\left(\dot{\boldsymbol{\Phi}}\, \boldsymbol{\Phi}^{-1}\right) = W \operatorname{tr} \mathbf{A} = 0, \tag{21}$$

and $\det \mathbf{M} = W(T) = W(0) = 1$. In fact, the monodromy matrix satisfies a much stronger property: it is symplectic. This means that

$$\mathbf{M}^\top \mathbf{J} \mathbf{M} = \mathbf{J}, \quad \mathbf{J} = \begin{bmatrix} 0_{N \times N} & I_{N \times N} \\ -I_{N \times N} & 0_{N \times N} \end{bmatrix}. \tag{22}$$

To prove it, note that $\mathbf{A}$ is a Hamiltonian matrix, i.e.,

$$\mathbf{A}^\top \mathbf{J} + \mathbf{J} \mathbf{A} = \mathbf{0}, \tag{23}$$

and expand

$$\overbrace{\boldsymbol{\Phi}^\top \mathbf{J} \boldsymbol{\Phi}}^{\cdot} = \dot{\boldsymbol{\Phi}}^\top \mathbf{J} \boldsymbol{\Phi} + \boldsymbol{\Phi}^\top \mathbf{J} \dot{\boldsymbol{\Phi}} = \boldsymbol{\Phi}^\top (\mathbf{A}^\top \mathbf{J} + \mathbf{J} \mathbf{A}) \boldsymbol{\Phi} = \mathbf{0}. \tag{24}$$

Consequently, $\boldsymbol{\Phi}^\top \mathbf{J} \boldsymbol{\Phi}$ is a constant matrix and is equal to its initial value of $\mathbf{J}$, and so is $\mathbf{M}^\top \mathbf{J} \mathbf{M}$.

## 3.2 Comments on numerical integration

If the modulated chain is stable, then all eigenvalues $\lambda$ of $\mathbf{M}$ lie inside the unit circle: $|\lambda| \leq 1$. Symplecticity further implies that the eigenvalues of $\mathbf{M}$ come in pairs of $\lambda$ and $1/\lambda$. Accordingly, if the chain is stable, then all eigenvalues lie exactly on the unit circle: $|\lambda| = 1$.

This is of numerical significance since the condition for stability is somewhat fragile. Indeed, any numerical amplification factors introduced by an explicit discretization scheme will lead to false instabilities. Thus, using explicit solvers, one needs to impose tight tolerances to keep such instabilities as mild as possible. This however necessitates the use of impractically small step sizes and can be computationally prohibitive. Implicit discretization schemes provide an alternative: they introduce a numerical damping making them more stable even for larger step sizes. Their drawback is that at each time step, a system of equations needs be solved. Even then, in the present context, implicit solvers appear to be more efficient than explicit solvers and particularly so for slow modulations, i.e., when computing $\mathbf{M}$ requires larger integration times.

Of particular interest here, is a set of implicit solvers qualified as "symplectic". Such solvers will provide an exactly symplectic monodromy matrix $\mathbf{M}$, up to roundoff errors, and are therefore expected to be better predictors of stability. Perhaps the simplest, non-trivial, example of a symplectic integration method is the symplectic Euler method [4]:

$$\begin{cases} u^{(i+1)} & = u^{(i)} + hM^{-1}(t_{i+1})p^{(i)}, \\ p^{(i+1)} & = p^{(i)} - hK(t_{i+1})u^{(i+1)}, \end{cases} \tag{25}$$

where $h$ is the time step. In matrix form, the scheme is

$$\mathbf{\Phi}^{(i+1)} = \begin{bmatrix} I & hM^{-1}(t_{i+1}) \\ -hK(t_{i+1}) & I - h^2K(t_{i+1})M^{-1}(t_{i+1}) \end{bmatrix} \mathbf{\Phi}^{(i)} \equiv \mathbf{\Phi}(t_i, t_{i+1})\mathbf{\Phi}^{(i)}. \tag{26}$$

It is straightforward to verify that the transition matrices $\mathbf{\Phi}(t_i, t_{i+1})$ are symplectic, regardless of $h$, and that the resulting approximation $\mathbf{M} = \prod_{i=0}^{N-1} \mathbf{\Phi}(t_i, t_{i+1})$ is symplectic as well.

The symplectic Euler method is implicit, technically speaking, but can be implemented just as easily as an explicit method and there is no need to solve any set of equations to complete a step. For that reason, this method is sometimes qualified as "semi-implicit" or "semi-explicit". Its downside is that it is of first order. Semi-implicit symplectic integrators of higher orders do exist but only for systems whose Hamiltonians are separable, e.g., $H(t, u, p) = V(t, u) + T(p)$. For modulated chains, the Hamiltonian is

$$H(t, u, p) = \frac{1}{2}u^\top K(t)u + \frac{1}{2}p^\top M^{-1}(t)p, \tag{27}$$

and is separable for instance if $M$ is constant. In these cases, high-order semi-implicit symplectic integrators are of preference. For more general modulations with non-separable Hamiltonians, the last resort is to use fully implicit symplectic integrators.

### Code: An implicit symplectic Runge-Kutta integrator

The following code implements a basic implicit Runge-Kutta solver, with a given Butcher tableau. The method being (fully) implicit, a system of equations is solved at each time step using a fixed-point iteration. When the equations are linear, as in the present context, the system can be solved directly (e.g., with Gaussian Elimination) but the fixed-point method was deemed more convenient: it is straightforward to implement, does not require the assembly of the linear system, and naturally takes advantage of the sparsity of matrix $\mathbf{A}$. Note in particular how the function `sys` below computes the matrix product $Ku$ without forming $K$. The monodromy matrix can then be computed by calling `irk` over `sys` with the identity matrix as an initial condition and one modulation period as the final time.

For reference, the following Butcher tableau is that of a symplectic implicit Runge-Kutta method, namely the 6th order Gauss-Legendre method:

$$\text{GL6} = \left[\begin{array}{c|ccc} \frac{1}{2} - \frac{\sqrt{15}}{10} & \frac{5}{36} & \frac{2}{9} - \frac{\sqrt{15}}{15} & \frac{5}{36} - \frac{\sqrt{15}}{30} \\[4pt] \frac{1}{2} & \frac{5}{36} + \frac{\sqrt{15}}{24} & \frac{2}{9} & \frac{5}{36} - \frac{\sqrt{15}}{24} \\[4pt] \frac{1}{2} + \frac{\sqrt{15}}{10} & \frac{5}{36} + \frac{\sqrt{15}}{30} & \frac{2}{9} + \frac{\sqrt{15}}{15} & \frac{5}{36} \\[4pt] \hline & \frac{5}{18} & \frac{4}{9} & \frac{5}{18} \end{array}\right]. \tag{28}$$

For a refresher on Runge-Kutta methods including `GL6`, see [5]. Here, with `GL6`, satisfactory results are obtained with 7 fixed-point iterations per step and 20 to 40 steps per characteristic time, a characteristic time being the period of the modulation or of an oscillation, whichever is smallest. Note again that when mass is time-independent, the more efficient semi-implicit methods should be preferred; see, e.g., [6].

```python
def irk(sys, ic, tf, steps, fp, butcher, *args):
    # solves dy/dt = sys(y, t, *args), y(0) = ic on [0, tf]
    # uses Implicit Runge-Kutta with tableau = butcher
    # ic = 2D-array initial conditions, tf = final time
    # steps = number of steps, fp = fixed-point iterations / step

    # step size
    h = tf / steps

    # initial conditions
    t = 0
    y = ic
    Y = [ic]

    # RK stages and slopes
    stages = np.shape(butcher)[0] - 1
    k = np.zeros((*np.shape(y), stages), y.dtype)

    # implicit Runge-Kutta
    for _ in range(steps):
        # solve iteratively for slopes
        for _ in range(fp):
            for s in np.arange(stages):
                k[..., s] = sys(y + h * k @ butcher[s, 1:],
                                t + butcher[s, 0] * h, *args)
        # update rule
        y = y + h * k @ butcher[stages, 1:]
        Y += [y]
        t += h

    return Y

def sys(phi, t, N, Z, nu, k, m):
    # returns \dot{\phi}=A(t)*\phi for a modulated chain
    # with parameters N, Z, nu, k, m
    # k, m = [average, amplitude]
```

```python
# initialize chain(t)
x = 2 * np.pi / Z * np.arange(N)
k = (k[0] + k[1] * np.cos(x - nu * t))[:, None]
m = (m[0] + m[1] * np.cos(x - nu * t))[:, None]

# compute K*u
U = phi[0:N, :]
T = k*(np.vstack([U[1:N,:],U[0,:]])-U)
KU = np.vstack([T[-1,:],T[0:N-1,:]])-T

return np.vstack((phi[N:2 * N, :] / m, -KU))
```

### 3.3 Reduced monodromy matrix

The fact that the modulation is progressive, or even space-dependent or not, was irrelevant to the computation of the monodromy matrix thus far. This is unsatisfactory and there are ways in which one can take advantage of the specific shape of the modulation. Indeed, note that even though the modulation is $T$-periodic, it suffices to wait a time period of $\tau = \xi/v = T/Z$ for the modulation profile to repeat, shifted through one lattice spacing. For instance, assuming the sequence of spring constants at $t = 0$ is

$$\underbrace{k_1, \ldots, k_{Z-1}, k_Z}_{\text{cell 1}}, \ldots, \underbrace{k_1, \ldots, k_{Z-1}, k_Z}_{\text{cell c}}, \tag{29}$$

then, at $t = \tau$, the sequence becomes

$$\underbrace{k_Z, k_1, \ldots, k_{Z-1}}_{\text{cell 1}}, \ldots, \underbrace{k_Z, k_1, \ldots, k_{Z-1}}_{\text{cell c}}. \tag{30}$$

The same applies to the sequence of masses. In matrix form, the progressiveness of the modulation translates into

$$M(t + \tau) = SM(t)S^\top, \quad K(t + \tau) = SK(t)S^\top, \quad S = \begin{bmatrix} 0 & & & 1 \\ 1 & \ddots & & \\ & \ddots & \ddots & \\ & & 1 & 0 \end{bmatrix}, \tag{31}$$

where $S$ is an $N \times N$ permutation matrix. Consequently,

$$\mathbf{A}(t + \tau) = \mathbf{S}\mathbf{A}(t)\mathbf{S}^\top, \quad \mathbf{S} = \begin{bmatrix} S & 0_{N \times N} \\ 0_{N \times N} & S \end{bmatrix}. \tag{32}$$

The above relationship implies that the long-term behavior of the system can be deduced from its behavior over a single reduced period $T/Z$, up to some spatial shifts.

This motivates the introduction of a reduced monodromy matrix $\mathbf{m} = \mathbf{\Phi}(T/Z)$. Then,

$$\mathbf{M} = \mathbf{S}^Z (\mathbf{S}^\top \mathbf{m})^Z. \tag{33}$$

To prove it, note that $\mathbf{m}$ propagates a solution from time $t = 0$ to time $t = T/Z$; then, $\mathbf{SmS}^\top$ further propagates the solution from time $t = T/Z$ to time $2T/Z$; $\mathbf{S}^2 \mathbf{mS}^{2\top}$ takes it from $2T/Z$ to $3T/Z$ and so on and so forth all the way up to time $t = T$. Composing these transformations yields the above factorization.

The fact that the monodromy matrix $\mathbf{M}$ admits such a factorization as a power of $\mathbf{m}$ has important implications, both practical and conceptual. For now, simply note that the factorization permits to reduce numerical integration times by a factor $Z$. Specifically, the reduced monodromy $\mathbf{m}$ can be computed by calling `irk` with a final time of $\tau$.

### 3.4 Dispersion diagrams: Theory

The eigenmodes $\varphi(t)$ of the modulated chain have initial conditions $\varphi(0)$ that are eigenvectors of the monodromy matrix, namely

$$\mathbf{M}\varphi(0) = \Lambda\varphi(0), \tag{34}$$

or equivalently, $\varphi(T) = \Lambda\varphi(0)$. This relationship holds across any interval of width $T$ so that $\varphi(t + T) = \Lambda\varphi(t)$ for all $t$. Indeed,

$$\varphi(t + T) = \mathbf{\Phi}(t)\mathbf{M}\varphi(0) = \Lambda\mathbf{\Phi}(t)\varphi(0) = \Lambda\varphi(t). \tag{35}$$

Consider now any plane wave component of $\varphi(t)$ of the form $\exp(ßqn)\exp(-ß\omega t)$, the above relationship then yields

$$\Lambda = \exp(-ß\omega T). \tag{36}$$

Hence, the eigenvalue $\Lambda$, also known as a Floquet multiplier, allows to determine an eigenfrequency $\omega$ modulo $2\pi/T$. Plotting the eigenfrequencies $\omega$ versus the wavenumbers $q$ of all eigenvectors $\varphi(0)$ produces a dispersion diagram.

Unfortunately, such a dispersion diagram would be highly ambiguous. The main drawback is that the diagram attributes the same $\omega$ to all of the Fourier components of a given eigenmode $\varphi(t)$. Indeed, note how in (36) eigenfrequency $\omega$ depends on the multiplier $\Lambda$ but not on the wavenumber $q$; not to mention that $\omega$ would only be known modulo $2\pi/T$. With the reduced monodromy matrix, one can do better, much better.

So let $\varphi(0)$ instead be an eigenvector of $\mathbf{S}^\top \mathbf{m}$, namely

$$\mathbf{S}^\top \mathbf{m}\varphi(0) = \lambda\varphi(0), \tag{37}$$

or equivalently, $\varphi(\tau) = \lambda \mathbf{S} \varphi(0)$. Here too, multiplier $\lambda$ applies across any interval of width $\tau$ since, by the same logic as before,

$$\varphi(t + \tau) = \mathbf{S}\boldsymbol{\Phi}(t)\mathbf{S}^\top \mathbf{m}\varphi(0) = \lambda \mathbf{S}\boldsymbol{\Phi}(t)\varphi(0) = \lambda \mathbf{S}\varphi(t). \tag{38}$$

In terms of nodal displacements and momenta, the above relationship reads

$$u_n(t + \tau) = \lambda u_{n-1}(t), \quad p_n(t + \tau) = \lambda p_{n-1}(t). \tag{39}$$

As a consequence, the eigenvectors of $\mathbf{S}^\top \mathbf{m}$ are initial conditions for waves that propagate progressively, much like the modulation, with profiles given by

$$u_n(t) = \lambda^n u_0(t - n\tau), \quad p_n(t) = \lambda^n p_0(t - n\tau). \tag{40}$$

Moreover, these initial conditions can be chosen to be simultaneous[4] eigenvectors of both $\mathbf{S}^\top \mathbf{m}$ and $\mathbf{S}^Z$. Hence, let $\varphi(0)$ further satisfy

$$\mathbf{S}^Z \varphi(0) = Q\varphi(0) \tag{41}$$

meaning

$$u_{n+Z}(t) = Qu_n(t), \quad p_{n+Z}(t) = Qp_n(t). \tag{42}$$

Since $u_{n+N} = u_n$, per the imposed periodic boundary conditions, it comes that $Q^{\mathtt{c}} = 1$, where $\mathtt{c}$ is the number of unit cells in the chain. Thus, there exists a $q = 2\pi\ell/N$ such that $Q = \exp(\beta q Z)$ and $0 \leq \ell < \mathtt{c}$. Now define the eigenfrequency $\omega$ to be such that

$$\exp(-\beta\omega\tau) = \lambda\exp(-\beta q) \tag{43}$$

and let

$$\tilde{u}_n(t) = u_n(t)\exp(-\beta q n)\exp(\beta\omega t). \tag{44}$$

Then, equations (40) and (42) together imply that

$$\tilde{u}_n(t) = \tilde{u}_0(t - n\tau) \tag{45}$$

and that $\tilde{u}_0$ is $T$-periodic in $t$ and $Z$-periodic in $n$.

To summarize: the simultaneous eigenvectors $\varphi(0)$ of $\mathbf{S}^T \mathbf{m}$ and $\mathbf{S}^Z$ propagate like Floquet-Bloch waves with

$$u_n(t) = \tilde{u}_0(t - n\tau)\exp(\beta q n)\exp(-\beta\omega t). \tag{46}$$

Finally, to produce a dispersion diagram, further expand $\tilde{u}_0$ into a Fourier series

$$\tilde{u}_0(t - n\tau) = \sum_j U_0^j \exp[\beta j \nu(n\tau - t)]. \tag{47}$$

---

[4] $\mathbf{S}^Z$ shifts through $Z$ lattice positions and thus leaves the chain as it is. In matrix form, this invariance property implies $\mathbf{S}^Z \mathbf{A} = \mathbf{A}\mathbf{S}^Z$ and subsequently $\mathbf{S}^Z \mathbf{m} = \mathbf{m}\mathbf{S}^Z$ and $\mathbf{S}^Z \mathbf{S}^\top \mathbf{m} = \mathbf{S}^\top \mathbf{m}\mathbf{S}^Z$. Commuting matrices admit simultaneous eigenvectors.

Then, the dispersion diagram is the locus of modes $(q + j\nu\tau, \omega + j\nu)$ such that $U_0^j$ is non-zero.

### 3.5 Dispersion diagrams: Numerical considerations

As argued above, producing a non-ambiguous dispersion diagram requires solutions be computed over one time period $[0, T]$, with initial conditions that are eigenvectors of both $\mathbf{S}^\top \mathbf{m}$ and $\mathbf{S}^Z$. There are two main ways in which this computation can be made more efficient.

On one hand, it is enough to find solutions over $[0, \tau]$ with initial conditions that are eigenvectors of $\mathbf{S}^\top \mathbf{m}$ then use the property $\varphi(t + \tau) = \lambda \mathbf{S} \varphi(t)$ to extend them to $[0, T]$.

On the other hand, instead of computing the $2N \times 2N$ shifted reduced monodromy matrix $\mathbf{S}^\top \mathbf{m}$ of the whole chain, one can impose the boundary conditions $u_Z = Q u_0$, $p_Z = Q p_0$, and compute the $2Z \times 2Z$ shifted reduced monodromy of a single unit cell. The system of ODEs written for a single unit cell has the same structure as (19) with its mass, constitutive and compatibility matrices given by

$$M = \begin{bmatrix} m_1 & & \\ & \ddots & \\ & & m_Z \end{bmatrix}, \quad C = \begin{bmatrix} k_1 & & \\ & \ddots & \\ & & k_Z \end{bmatrix}, \quad A = \begin{bmatrix} -1 & 1 & & \\ & \ddots & \ddots & \\ & & & 1 \\ Q & & & -1 \end{bmatrix}. \quad (48)$$

The expression of matrix $A$ embodies the boundary conditions. As a result, $A$ is now complex-valued and the stiffness matrix $K = A^\dagger C A$ refers to its conjugate transpose $A^\dagger$ rather than to its transpose $A^\top$. Similarly, the permutation matrices become

$$S = \begin{bmatrix} 0 & & & Q^* \\ 1 & \ddots & & \\ & \ddots & \ddots & \\ & & 1 & 0 \end{bmatrix}, \quad \mathbf{S} = \begin{bmatrix} S & 0_{Z \times Z} \\ 0_{Z \times Z} & S \end{bmatrix}, \quad (49)$$

and the shifted reduced monodromy is $\mathbf{S}^\dagger \mathbf{m}$ with $\mathbf{S}^\dagger$ being the conjugate transpose of $\mathbf{S}$. The system of ODEs written for a single unit cell must be solved for all $c$ possible values of $Q$. Thus, the gain in computational complexity is equal to the ratio of $O(N^\alpha)$ to $O(\mathtt{c}Z^\alpha)$, namely $O(\mathtt{c}^{\alpha-1})$, assuming $O(N^\alpha)$ for some $\alpha$ is the computational complexity of finding the monodromy matrix and its eigenvectors for a chain of length $N$.

### Code: The dispersion diagram of a modulated chain

The system of ODEs for a single unit cell can be solved using the implicit solver `irk` called over a modified `sys` function implemented below. Then, the derivations above are coded into a function that computes dispersion diagrams.

```python
def sys(phi, t, Q, Z, nu, k, m):
    # returns \dot{\phi}=A(t)*\phi for a single unit cell
    # of a modulated chain with parameters Z, nu, dk, dm
    # and Floquet multiplier Q

    # initialize chain(t)
    x = 2 * np.pi / Z * np.arange(Z)
    k = (k[0] + k[1] * np.cos(x - nu * t))[:, None]
    m = (m[0] + m[1] * np.cos(x - nu * t))[:, None]

    # compute K*u
    U = phi[0:Z, :]
    T = k * (np.vstack([U[1:Z, :], Q * U[0, :]]) - U)
    KU = np.vstack([T[-1, :] / Q, T[0:Z - 1, :]]) - T

    return np.vstack((phi[Z:2 * Z, :] / m, -KU))

def dispersion(Z, c, nu, k, m, substeps=20, fp=7):
    # returns the dispersion diagram
    # (Wnb: Wavenumbers, Freq: Frequencies, Amp: Amplitudes)
    # of a modulated chain (Z, c, nu, k, m)
    # substeps and fp are numerical integration parameters

    # characteristic times
    T = 2 * np.pi / nu
    tau = T / Z
    T_nat = 2 * np.pi * np.sqrt((m[0] - m[1]) / (k[0] + k[1]))

    # number of steps for numerical integration
    steps = int(tau / min(tau, T_nat)) * substeps

    # wavenumbers
    q = np.arange(0, 2 * np.pi / Z, 2 * np.pi / Z / c)

    # initialize dispersion diagram
    Wnb  = np.zeros((Z * steps * c, 2 * Z))
    Freq = np.zeros((Z * steps * c, 2 * Z))
    Amp  = np.zeros((Z * steps * c, 2 * Z))

    for ell in range(c):
        # integrate over [0, tau] + boundary condition: u_Z = Q u_0
        Q = np.exp(1j * q[ell] * Z)
        Phi = irk(sys, np.eye(2 * Z, dtype=np.cdouble), tau, steps,
                  fp, GL6, Q, Z, nu, k, m)

        # reduced monodromy
        mon = Phi.pop()

        # shifted reduced monodromy
        newOrder = np.roll(np.arange(2 * Z).reshape(2, Z), -1,
                           axis=1).flatten()
        stm = mon[newOrder, :]
        stm[[Z - 1, -1], :] = stm[[Z - 1, -1], :] * Q
```

```python
    # eigen-values, vectors and frequencies
    ev, EV0 = np.linalg.eig(stm)
    omega = np.angle(np.exp(1j * q[ell]) / ev) / tau

    # propagate Floquet modes over [0, tau[
    EVt = Phi @ EV0

    # extend propagation of mass 0 to one full period [0, T[
    u0 = EVt[:, 0, :]
    for j in range(1, Z):
        u0 = np.concatenate((u0, EVt[:, Z - j, :] * ev**j / Q))

    # extend time to [0, T[
    t = np.linspace(0, T, Z * steps + 1)[:-1]

    # extract periodic displacement profile
    utilde0 = u0 * np.exp(1j * omega * t[:, None])

    # compute Fourier components, normalized
    U0 = np.abs(np.fft.fft(utilde0, axis=0))
    U0 = U0 / np.max(U0, axis=0)

    # frequencies and wavenumbers
    cutoff = 2 * np.pi * Z * steps / T
    J = np.arange(Z * steps, 0, -1)[:, None]
    freq = np.mod(J * nu + omega, cutoff)
    wnb = np.mod(2 * np.pi * J / Z + q[ell], 2 * np.pi)

    # re-center on 1st Brillouin zone
    freq[freq > cutoff / 2] = freq[freq > cutoff / 2] - cutoff
    wnb[wnb > np.pi] = wnb[wnb > np.pi] - 2 * np.pi

    # store diagram
    Freq[ell * Z * steps:(ell + 1) * Z * steps, :] = freq
    Wnb[ell * Z * steps:(ell + 1) * Z * steps, :] = wnb
    Amp[ell * Z * steps:(ell + 1) * Z * steps, :] = U0

return Wnb, Freq, Amp
```

## 3.6 What non-reciprocity looks like

With the functions written thus far, it is possible to compute accurate dispersion diagrams and simulate wave propagation over long periods of time with relative comfort. Four cases are illustrated on Fig. 4. The homogeneous chain exhibits its signature sine-shaped dispersion curve. Inhomogeneities break that curve into several branches, four in this case. Indeed, the exemplified chain is periodic with 4 masses per unit cell. Either way, the dispersion diagram is even, symmetric with respect to the frequency axis: waves going left and right propagate at the same frequencies with identical dispersive properties. Conversely, the space-time diagrams appear to be symmetric enough: rays emanated left and right are equally likely and propagate at equal speeds and amplitudes. Enter the modulation: the symmetry of the dispersion diagram breaks down. Now, depending on the sign of the wavenumber, harmonic components will travel at higher or lower frequencies with higher or lower group and phase velocities. On the space-time diagrams, it is somewhat clear that rays emanated to the right are slower, more powerful and more frequent than rays emanated to the left. For faster modulations, the effects become even more dramatic. This left-right bias in wave propagation observed on Fig. 4 is loosely[5] referred to as a "non-reciprocity".

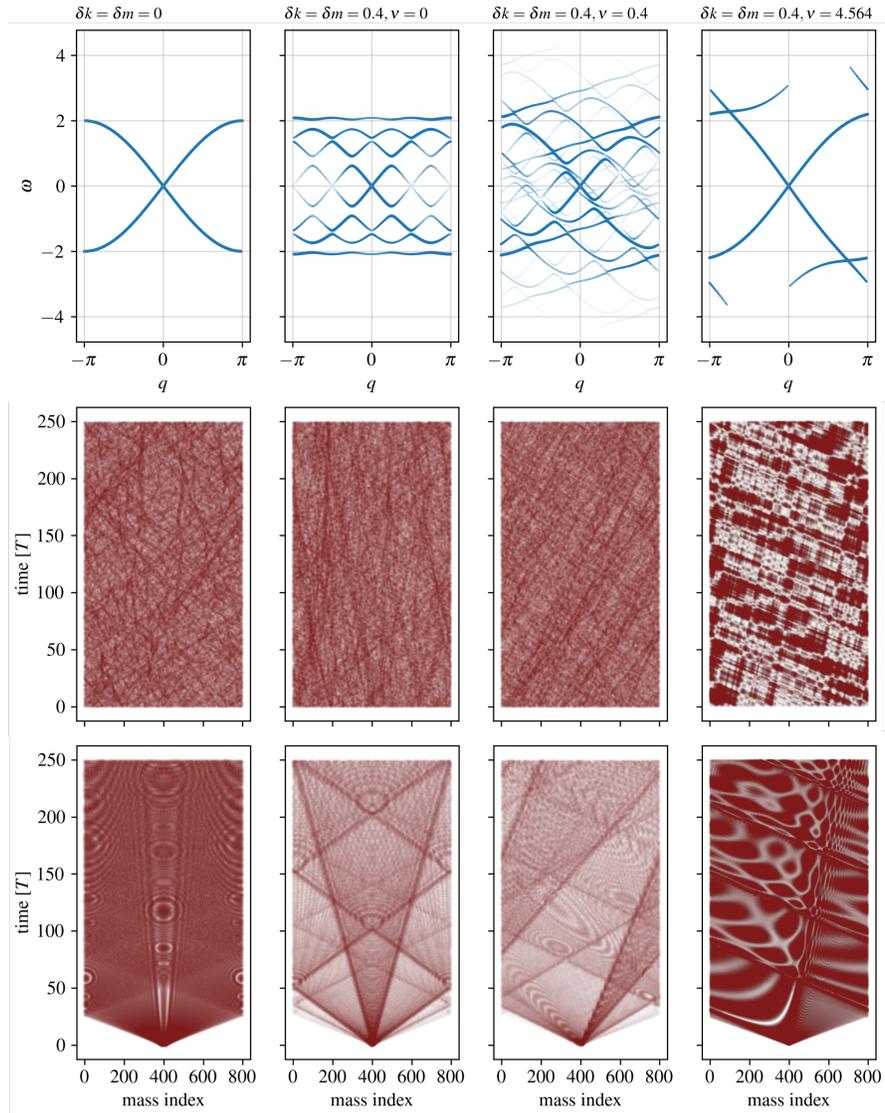

**Fig. 4** Dispersion and space-time diagrams in non-modulated and modulated chains. Left to right, by columns: a homogeneous non-modulated chain, a periodic non-modulated chain, a modulated chain, a modulated chain on the brink of instability. Top to bottom, by rows: dispersion diagrams, space-time diagrams of field intensity for a random normally-distributed initial condition, space-time diagrams for a displacement Dirac initial condition. Point size and transparency are correlated with amplitude. Other numerical parameters are: $Z = 4$, $c = 30$, $k_o = m_o = 1$. Time is measured in units of $T = 2\pi/0.4$.

## 4 Progressively modulated continua

One can easily be persuaded that progressive modulations should cause non-reciprocity: waves going with or against the modulation are scattered differently by a Doppler-like effect. What non-reciprocity will exactly look like, on the other hand, is harder to predict as seen on Fig. 4. Progressively modulated continua constitute a remarkable exception: hereafter, we deduce, in closed form, the influence of a progressive periodic modulation on wave propagation in a 1D, originally non-dispersive, continuum. The main result is a transformation mapping waves in the modulated medium to waves in a fictitious non-modulated medium. Last, non-reciprocal effects for waves of low frequencies are investigated and brief specialized comments pertaining to "Willis coupling" are presented.

### 4.1 The continuum limit

The reduced monodromy matrix leveraged the progressiveness of the modulation and permitted to shorten integration times from $T$ to $\tau = T/Z$ where $Z$ is the size of a unit cell. In the continuum limit, $Z$ approaches infinity and the integration time $\tau$ approaches 0. In other words, for a progressively modulated continuum, the reduced monodromy, and therefore the dispersion diagram, can be computed without the need for time integration at all. As a matter of fact, for $\tau \to 0$,

$$\varphi(\tau) = \varphi(0) + \tau\dot{\varphi}(0) = \varphi(0) + \tau\mathbf{A}\varphi(0), \tag{50}$$

meaning that

$$\mathbf{m} = \mathbf{I} + \tau\mathbf{A} \tag{51}$$

to first order in $\tau$. Thus, computing the dispersion diagram amounts to solving the eigenvalue problem

---

[5] For a technical statement of the principle of reciprocity, see [2].

$$\varphi(0) + \tau \mathbf{A}\varphi(0) = \lambda \mathbf{S}\varphi(0) \tag{52}$$

and involves no integration in time.

Keep in mind that as $\tau$ approaches 0, the reduced monodromy matrix grows in size proportionately to $Z$ and gains new eigenmodes with increasingly high frequencies. Here, it is understood that the frequencies, as well as the wavenumbers, of interest do not grow unbounded with $Z$. In particular, as $\tau$ approaches 0, it is assumed that $\lambda$ approaches 1 and $\mathbf{S}\varphi(0)$ approaches $\varphi(0)$. Formally, both $\lambda$ and $\mathbf{S}$ expand into

$$\lambda = 1 - \text{\ss}\tau\Omega, \quad \mathbf{S} = \mathbf{I} + a\boldsymbol{\Delta}^{\top}, \tag{53}$$

where $\Omega$ is a finite frequency that is to be determined and interpreted; $\boldsymbol{\Delta} = \frac{\mathbf{S}^{\top} - \mathbf{I}}{a}$ is a discretized space differential operator; and, $a$ is a small distance separating two consecutive masses.

Then, to leading order, the eigenvalue problem transforms into

$$\mathbf{A}\varphi(0) = -\text{\ss}\Omega\varphi(0) + c\boldsymbol{\Delta}^{\top}\varphi(0), \tag{54}$$

where $c = a/\tau$ is the modulation speed. It is worth noting that the dispersion diagram of the non-modulated chain can be determined from the standard eigenvalue problem

$$\mathbf{A}\varphi(0) = -\text{\ss}\Omega\varphi(0). \tag{55}$$

Thus, overall, the effect of the modulation is to introduce a drift term $c\boldsymbol{\Delta}^{\top}$ in the form of a space derivative proportional to the modulation speed. Most importantly, when the modulation changes direction or, equivalently, when a propagated wave is incident in the opposite direction, the drift $c\boldsymbol{\Delta}^{\top}$ changes sign. This drift is accordingly at the origin of the left-right bias and of the emerging non-reciprocity.

The question investigated in what follows is: how does the progressive modulation modify solutions of (55) into solutions of (54)?

## 4.2 A change in perspective

Thus far, it was natural to consider the equations of motion as ODEs in time. In the continuum limit however, time integration is no longer necessary. It is then of interest to change perspective and to rewrite the equations of motion as ODEs in space. Indeed, solving equation (54) amounts to solving a set of two ODEs

$$\begin{aligned} p/m &= -\text{\ss}\Omega u - cu', \\ a^2(ku')' &= -\text{\ss}\Omega p - cp'. \end{aligned} \tag{56}$$

To see why, recall that $K = A^{\top}CA$ and note that $A = S^{\top} - I = a\Delta$ can be identified with $a\partial/\partial x$ and $A^{\top}$ can be identified with its adjoint $-a\partial/\partial x$, with $\partial/\partial x \equiv \cdot'$ being the space derivative. The unknown fields $u = u(x)$ and $p = p(x)$ are the continuous limits of $u_n(0)$ and $p_n(0)$ where $na \to x$; together, they define initial conditions for the eigenmodes of the progressively modulated continuum. Furthermore, it is convenient to redefine $m$ in terms of a mass density $\rho = m/a$; to replace linear momentum $p$ with a linear momentum density $p^\infty = p/a$; and, to redefine spring constant $k$ in terms of a "string tension" $k^\infty = ak$. In what follows, notations are abused and the superscript $\infty$ is dropped.

In conclusion, the ODEs of interest are

$$
\begin{aligned}
p/\rho &= -\text{ß}\Omega u - cu', \\
(ku')' &= -\text{ß}\Omega p - cp',
\end{aligned}
\tag{57}
$$

where $\rho = \rho(x)$ and $k = k(x)$ are periodic in $x$ and describe the initial state of the modulation. Once the initial states $u$ and $p$ are determined, they can be propagated in time according to (40), namely

$$
u(x,t) = u(x - ct)\exp(-\text{ß}\Omega t), \quad p(x,t) = p(x - ct)\exp(-\text{ß}\Omega t).
\tag{58}
$$

### 4.3 The main result

It is convenient to write the system of ODEs as a pair of first order ODEs for the particle velocity $v = p/\rho$ and stress $\sigma = ku$,

$$
\begin{aligned}
v' + c(\sigma/k)' &= -\text{ß}\Omega\sigma/k, \\
\sigma' + c(\rho v)' &= -\text{ß}\Omega\rho v.
\end{aligned}
\tag{59}
$$

This suggests the change of variables

$$
\psi = \begin{bmatrix} v + c\sigma/k \\ \sigma + c\rho v \end{bmatrix}
\tag{60}
$$

for which

$$
\psi' = -\text{ß}\Omega \begin{bmatrix} 0 & 1/k \\ \rho & 0 \end{bmatrix} \begin{bmatrix} v \\ \sigma \end{bmatrix} = \frac{\text{ß}\Omega c}{s^2 - c^2}\psi - \frac{\text{ß}\Omega s^2}{s^2 - c^2} \begin{bmatrix} 0 & 1/k \\ \rho & 0 \end{bmatrix} \psi
\tag{61}
$$

with $s(x) = \sqrt{k(x)/\rho(x)}$ being the speed of sound at position $x$. This last form invites a further change of variables

$$
\begin{bmatrix} V \\ \Sigma \end{bmatrix} = \psi \exp\left( -\text{ß}\Omega \int_0^x \frac{c\,\mathrm{d}x}{s^2 - c^2} \right),
\tag{62}
$$

which leads to an even simpler set of ODEs, namely

$$V' = -\beta\Omega\,\Sigma/\bar{k}, \qquad \bar{k} \equiv k\,\frac{s^2 - c^2}{s^2},$$

$$\Sigma' = -\beta\Omega\,\bar{\rho}V, \qquad \bar{\rho} \equiv \rho\,\frac{s^2}{s^2 - c^2}. \tag{63}$$

The major advantage of the lastly adopted change of variables is that it turns the modulation speed $c$ into a mere numerical coefficient whose presence or absence does not alter the form of the system of ODEs. In particular, it permits to state the main result of the present section: *the eigenmodes of the progressively modulated medium of properties $(k, \rho)$ are in a one-to-one correspondence with the eigenmodes of a fictitious non-modulated medium of properties $(\bar{k}, \bar{\rho})$*. Specifically, given an eigenmode of the fictitious non-modulated medium, the eigenmode (58) of the modulated medium is

$$v(x, t) = v(x_t)\exp(-\beta\Omega t)$$

$$= \left(V(x_t) - c\,\frac{\Sigma(x_t)}{k(x_t)}\right)\frac{s^2(x_t)}{s^2(x_t) - c^2}\exp\left[-\beta\Omega\left(t - \int_0^{x_t}\frac{c\,\mathrm{d}x}{s^2 - c^2}\right)\right] \tag{64}$$

with $x_t \equiv x - ct$. In particular, if the eigenmode of the fictitious medium has frequency $\Omega$ and wavenumber $Q$, then the eigenmode of the modulated medium has frequency $\omega$ and wavenumber $q$ given by

$$\omega = \Omega + \Omega\left\langle\frac{c^2}{s^2 - c^2}\right\rangle + cQ, \quad q = Q + \Omega\left\langle\frac{c}{s^2 - c^2}\right\rangle \tag{65}$$

where $\langle\rangle$ denotes the average value over a unit cell.

The above result is valid as long as the change of variables is one-to-one, i.e., as long as the modulation speed is different from the speed of sound at all positions $x$. Note that the fictitious material will exhibit negative values for $\bar{k}$ and $\bar{\rho}$ wherever the modulation speed exceeds the speed of sound. Note also that $\bar{k}\bar{\rho} = k\rho$, implying $\bar{z} = z$ where $z = \rho s$, $\bar{z} = \bar{\rho}\bar{s}$ are the impedances and $\bar{s} = \sqrt{\bar{k}/\bar{\rho}}$ is the associated modulated wave speed. The fictitious medium therefore has the same impedance as the non-modulated medium and differs only in terms of the reduced local sound speed $\bar{s} = s - c^2/s$.

Finally, it is noteworthy that the above correspondence is not restricted to real values of $(Q, \Omega)$ and $(q, \omega)$. In particular, it would be inaccurate in principle to deduce that the modulated medium is stable as soon as the fictitious medium is stable. Indeed, the fictitious medium could be stable but still allow for modes where both $Q$ and $\Omega$ are complex in such a way that $q$ is real but $\omega$ is complex. In these cases, the modulated medium would be unstable.

### 4.4 What non-reciprocity looks like

The dispersion diagram of a progressively modulated continuum is exemplified on Fig. 5. For reference, the diagrams of the original non-modulated continuum and of the corresponding fictitious continuum are shown as well. As the modulation speed $c$ increases and approaches the speed of sound $s$, the speed of sound $\bar{s}$ in the fictitious medium is reduced further causing the dispersion diagram to display lower group velocity (which is always bounded above by the sound speed), with the effect that the diagram is compressed in the vertical direction. Following that action, the dispersion diagram of the modulated continuum can be obtained by the shearing transformation (65). This invariably breaks the parity of the diagram and introduces a non-reciprocal left-right bias, which is evident from the group velocity relation

$$\frac{\mathrm{d}\omega}{\mathrm{d}q} = c + \frac{\mathrm{d}\Omega}{\mathrm{d}Q} \left/ \left( 1 + \left\langle \frac{c}{s^2 - c^2} \right\rangle \frac{\mathrm{d}\Omega}{\mathrm{d}Q} \right) \right. \tag{66}$$

that breaks the positive/negative symmetry of $\mathrm{d}\Omega/\mathrm{d}Q$ in favor of positive $\mathrm{d}\omega/\mathrm{d}q$. Note that the shown numerical results were obtained through direct time integration using the techniques deployed in section 3.

Most noteworthy is the fact that for a strong-enough and fast-enough modulation, the shear transformation can be so severe that the acoustic branches will exhibit two group velocities of the same sign. This suggests that low-frequency long-wavelength propagation will be confined to one direction and prohibited in the opposite direction. In other words, a localized source will emit two rays in the same direction, namely to the right in the exemplified case. Time-domain simulations shown on the same figure confirm the prediction. Hereafter, this non-reciprocal effect[6] is characterized more closely.

### 4.5 Acoustic group velocities and Willis coupling

At low frequencies and long wavelengths, waves in the fictitious medium propagate at an acoustic group velocity $\hat{s}$ given by

$$\hat{s}^2 = \frac{1}{\langle \bar{\rho} \rangle \left\langle 1/\bar{k} \right\rangle}, \tag{67}$$

that characterizes the slope of $\Omega$ in function of $Q$. In other words, $\Omega/Q = \pm\hat{s}$. Accordingly, the acoustic group velocity in the modulated medium is

$$s_\pm \equiv \frac{\omega}{q} = \frac{c \pm \hat{s} \left\langle \frac{s^2}{s^2 - c^2} \right\rangle}{1 \pm \hat{s} \left\langle \frac{c}{s^2 - c^2} \right\rangle} \tag{68}$$

and can take one of two values. In the absence of modulation ($c = 0$), these two values are $s_\pm = \pm \hat{s}$ and are equal and opposite. The modulation introduces a bias that, to leading order in $c$, reads

$$\begin{aligned} s_\pm &= \pm \hat{s}_o + \left( \langle \rho \rangle \langle 1/k \rangle - \langle \rho/k \rangle \right) \hat{s}_o^2 \, c \\ &= \pm \hat{s}_o - \left\langle (\rho - \langle \rho \rangle) \left( 1/k - \langle 1/k \rangle \right) \right\rangle \hat{s}_o^2 \, c, \end{aligned} \tag{69}$$

where $\hat{s}_o = 1/(\langle \rho \rangle \langle 1/k \rangle)$ is the effective speed of sound in the absence of modulation. Hence, when $\rho$ and $1/k$ are correlated (resp. anti-correlated) so that their covariance is positive (resp. negative), slow modulations will decelerate (resp. accelerate) co-propagated waves and accelerate (resp. decelerate) counter-propagated waves. Note that in the example above, $k = s^2 \rho$ where $s$ is constant implies that $\rho$ and $1/k$ are anti-correlated. In any case, $s_\pm$ maintain opposite signs.

---

[6] This effect is described in [7] in a particular case and is referred to as "coordinated wave propagation".

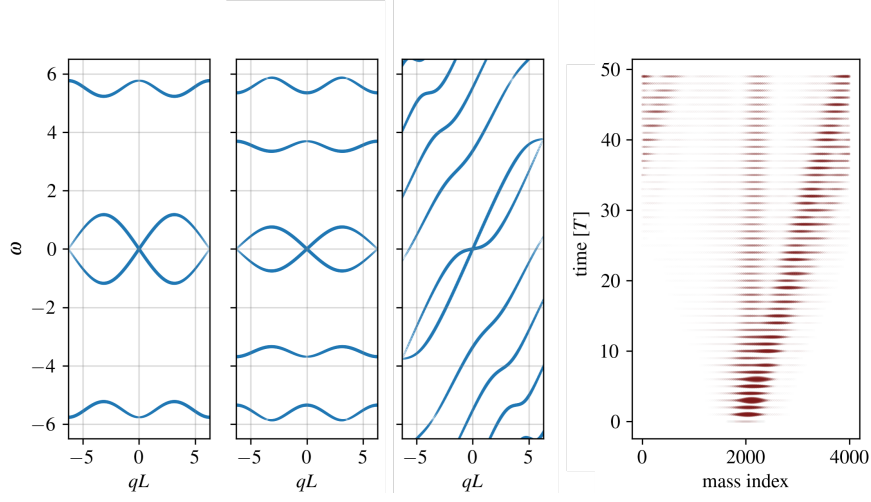

**Fig. 5** Wave propagation in a progressively modulated continuum. Left to right: dispersion diagrams of the original non-modulated medium, of the corresponding fictitious non-modulated medium, and of the progressively modulated medium, respectively; space-time diagram of field intensity for a long-wavelength Gaussian-shaped initial condition. Numerical parameters are: $\mathtt{c} = 100$, $k_o = m_o = 1$, $\delta k = \delta m = 0.95$, $c = 0.6$. Time is measured in units of $T = 2\pi L/c$. The number of masses per unit cell, $Z$, is infinite in theory; satisfactory results in reasonable times were obtained with $Z = 20, 20, 100, 40$ respectively.

By contrast, for faster modulations, and in particular, when $c$ approaches the speed of sound of one of the constitutive phases, it can be seen that $\hat{s}$ approaches 0 so that $s_+$ and $s_-$ both approach $c$. Thus, for some critical value of the modulation speed, between 0 and $\min s$, one of the acoustic group velocities flips sign. As a result, the modulated medium exhibits two co-propagating acoustic modes. Specifically, $s_-$ changes sign when $c$ first exceeds the smallest value satisfying

$$\left\langle \frac{s^2}{s^2 - c^2} \right\rangle^2 = c^2 \left\langle \frac{zs}{s^2 - c^2} \right\rangle \left\langle \frac{s}{z(s^2 - c^2)} \right\rangle. \tag{70}$$

For example, the critical modulation speed for constant sound speed is $c = s/\sqrt{\langle z \rangle \langle 1/z \rangle}$, which, as expected, lies below $s$ in value.

The left-right bias, be it extreme or not, implies that the usual non-dispersive wave equation governing the propagation of acoustic modes, namely

$$k^e u_{,xx} - \rho^e u_{,tt} = 0, \tag{71}$$

breaks down. Indeed, such equations are bound to produce $s_+ = -s_-$. Here, the low-frequency dispersion equation

$$(\omega - s_+ q)(\omega - s_- q) = 0, \tag{72}$$

i.e.,

$$\omega^2 - (s_+ + s_-)\omega q + s_- s_+ q^2 = 0, \tag{73}$$

where $s_+ + s_- \neq 0$ suggests that a more appropriate wave equation is

$$k^e u_{,xx} + 2W u_{,xt} - \rho^e u_{,tt} = 0. \tag{74}$$

The mixed derivative coefficient, $W$, is known as a Willis coupling and is directly responsible for the discrepancy between acoustic group velocities. In terms of an effective constitutive law, $k^e$ is interpreted as an effective elasticity modulus, $\rho^e$ is an effective mass density and $W$ is a coupling term between stress and velocity as well as between momentum and strain. For a detailed derivation of such an effective model and a generalization to three space dimensions, see [8].

# 5 Concluding remarks

The theoretical and numerical elements of this chapter should help to model and to computationally solve systems with space-time dependent material properties defined by a periodic progressive modulation. Starting from a single DOF system modeled by Mathieu's equation, and then moving on to a chain of a finite number of springs and masses governed by a system of ODEs, it became evident that the fundamental quantity of interest for determining stable solutions is the monodromy matrix defined by the unit time period. In the case of the finite chain with multiple unit cells and periodic end conditions, the simpler reduced monodromy matrix for the unit cell naturally falls out from the analysis.

The monodromy matrix concept is widespread in physics and engineering, particularly in systems with material properties that are spatially periodic and independent of time. In that case the monodromy matrix is the *matricant* evaluated over a unit spatial period, a quantity that can be evaluated using standard methods of linear algebra [9]. By contrast, computation of the monodromy matrix for systems with progressive space-time modulation requires integration of a set of ODEs over the unit time period. A specific code has been described for that purpose, using a high order Runge-Kutta scheme that is also symplectic, guaranteeing that the numerical solution conserves quantities analogous to a wave energy flux. Details on how to compute dispersion curves (band diagrams) are also provided with explicit code for that purpose. The behavior of a progressively modulated continuum can then be analyzed by taking the appropriate limits.

The breakdown of wave reciprocity has been demonstrated for systems with both finite and infinite DOFs through explicit computations showing preferential propagation in one direction. Wave reciprocity is a fundamental physical result in materials with linear time-independent properties, even in the presence of damping [2]. Systems that violate reciprocity are of interest for their ability to allow mechanical energy to travel in one direction but not its reverse. The methods and codes provide a template for modeling non-reciprocity and computing solutions in more complex and physically realistic systems.

## 6 Further reading

In preparing the material on Mathieu's equation and on the properties of the monodromy matrix, references [10,11] were particularly helpful. The reader interested in the general theory of time-periodic ODEs as well as in other solution strategies besides the monodromy matrix, including asymptotic analysis and Hill's determinants, is invited to consult them. For examples of band diagrams in modulated materials using Floquet–Bloch analysis and time-domain simulations, see [8,12–16]. It is generally observed that the spatial profile of the modulation dictates the band structure, while the degree of band shearing depends upon the modulation frequency.

For the origins of the Willis material model, see [17,18]. As noted above, Willis coupling introduces a first-order time derivative to the wave equation which breaks time-reversal symmetry and, consequently, reciprocity [7,19,20].

The models considered in this chapter correspond to *activated materials* which are materials with constitutive properties modulated in response to an external stimulus. Examples include a phononic crystal with bulk elastic moduli periodically modulated in space and time [15], a two-phase laminate with interfaces moving at the modulation velocity [8], and a metabeam with constant elastic moduli while the stiffness of locally resonant attachments is modulated [21]. A related Doppler-like effect occurs when the material properties are modulated in space and time by the action of a *pump wave* [22]; tuning the Doppler shift can achieve one-way Bragg reflection [23,24]. Activated non-reciprocity is also possible in media with multiple wave modes, where for instance an incident acoustic mode is directionally reflected into an optical mode [25]. Other modulated platforms for activated non-reciprocity include 1D piezoelectric structures with space-time varying electrical boundary conditions [26, 27], and media with modulated effective mass [28] or effective stiffness [29].

Modulation of material properties using *slow* pump waves can lead to significant accumulated Doppler-like phase shift over a full period of the modulation cycle. Effects such as band shearing and tilting can then be understood in terms of a *Berry phase*, a quantity originally introduced to explain quantal interference phenomena caused by slowly changing environmental parameters [30]. In particular, the band tilts by an amount proportional to the total *Berry's curvature* [8, 15], related to a quantized topological invariant known as the Chern number that is immune to small perturbations of the medium properties. A related topological phenomenon in time modulated materials is the occurrence of one-way edge modes circulating either clockwise or anticlockwise along the boundary [15, 31, 32].

The limit of *infinite* modulation speed corresponds to material properties that are time dependent and independent of spatial position. This results in a band-gap amplification effect, first noted in the context of electrical transmission lines with time-varying capacitance [33]. Recent studies have considered time-dependent acoustic [34] and elastic media [35], as well as space-time checkerboard patterns with novel wave effects [36]. For the effect of damping in stabilizing the exponential parametric growth, see [37].

## References

1. Brillouin, L.: Wave Propagation in Periodic Structures. Dover, New York (1953)
2. Nassar, H. *et al.*: Nonreciprocity in acoustic and elastic materials. Nature Rev. Mat. (2020) doi: 10.1038/s41578-020-0206-0
3. Mathieu, E.: Mémoire sur le mouvement vibratoire d'une membrane de forme elliptiques. J. Math. Pures Appl. **13**, 137–203 (1868)
4. Hairer, E., Lubich, C., Wanner, G.: Geometric Numerical Integration: Structure-Preserving Algorithms for Ordinary Differential Equations. Springer, Heidelberg (2006)
5. Iserles, A.: A First Course in the Numerical Analysis of Differential Equations. Cambridge University Press (1996)
6. Yoshida, H.: Construction of higher order symplectic integrators. Phys. Let. A (1990) doi: 10.1016/0375-9601(90)90092-3
7. Lurie, K. A.: An Introduction to the Mathematical Theory of Dynamic Materials. Springer, New York (2007)
8. Nassar, H., Xu, X. C., Norris, A. N., Huang, G. L.: Modulated phononic crystals: Non-reciprocal wave propagation and Willis materials. J. Mech. Phys. Solids, (2017) doi: 10.1016/j.jmps.2017.01.010
9. Pease, M.C.: Methods of Matrix Algebra. Academic Press, New York (1965)
10. Kovacic, I., Rand, R., Sah, S. M.: Mathieu's equation and its generalizations: overview of stability charts and their features. App. Mech. Rev. (2018) doi: 10.1115/1.4039144
11. Stoker, J. J.: Nonlinear vibrations in Mechanical and electrical systems. Interscience Publishers, New York (1950).
12. Vila, J., Pal, R. K., Ruzzene, M., Trainiti, G.: A Bloch-based procedure for dispersion analysis of lattices with periodic time-varying properties. J. Sound Vib. (2017) doi: 10.1016/j.jsv.2017.06.011

13. Wallen, S. P., Haberman, M. R.: Nonreciprocal wave phenomena in spring-mass chains with effective stiffness modulation induced by geometric nonlinearity. Phys. Rev. E (2019) doi: 10.1103/PhysRevE.99.013001

14. Goldsberry, B. M., Wallen, S. P., Haberman, M. R.: Non-reciprocal wave propagation in mechanically-modulated continuous elastic metamaterials. J. Acoust. Soc. Am. (2019) doi: 10.1121/1.5115019

15. Nassar, H., Chen, H., Norris, A., Huang, G.: Quantization of band tilting in modulated phononic crystals. Phys. Rev. B (2018) doi: 10.1103/PhysRevB.97.014305

16. Attarzadeh, M. A., Nouh, M.: Elastic wave propagation in moving phononic crystals and correlations with stationary spatiotemporally modulated systems. AIP Adv. (2018) doi: 10.1063/1.5042252

17. Willis, J. R. Variational principles for dynamic problems for inhomogeneous elastic media. Wave Motion, (1981) doi: https://doi.org/10.1016/0165-2125(81)90008-1

18. Willis, J. R. Dynamics of composites. In: Suquet, P. (ed.) Continuum Micromechanics, pp. 265–290. Springer-Verlag, New York (1997)

19. Lurie, K. A.: Effective properties of smart elastic laminates and the screening phenomenon. Int. J. Solids Struct. (1997) doi: 10.1016/S0020-7683(96)00105-9

20. Quan, L., Sounas, D. L., Alù, A.: Nonreciprocal Willis coupling in zero-index moving media. Phys. Rev. Lett. (2019) doi: 10.1103/PhysRevLett.123.064301

21. Nassar, H., Chen, H., Norris, A. N. Huang, G. L. Non-reciprocal flexural wave propagation in a modulated metabeam. Extr. Mech. Lett. (2017) doi: 10.1016/j.eml.2017.07.001

22. Slater, J. C.: Interaction of waves in crystals. Rev. Mod. Phys. (1958) doi: 10.1103/RevModPhys.30.197

23. Swinteck, N. *et al.*: Bulk elastic waves with unidirectional backscattering-immune topological states in a time-dependent superlattice. J. Appl. Phys. (2015) doi: 10.1063/1.4928619

24. Trainiti, G., Ruzzene, M.: Non-reciprocal elastic wave propagation in spatiotemporal periodic structures. New J. Phys. (2016) doi: 10.1088/1367-2630/18/8/083047

25. Nassar, H., Chen, H., Norris, A. N., Haberman, M. R., Huang, G. L.: Non-reciprocal wave propagation in modulated elastic metamaterials. Proc. R. Soc. A (2017) doi: 10.1098/rspa.2017.0188

26. Croënne, C., Vasseur, J. O., Matar, O. B., Hladky-Hennion, A.-C., Dubus, B.: Non-reciprocal behavior of one-dimensional piezoelectric structures with space-time modulated electrical boundary conditions. J. Appl. Phys. (2019) doi: 10.1063/1.5110869

27. Marconi, J. *et al.*: Experimental observation of non-reciprocal band-gaps in a space-time modulated beam using a shunted piezoelectric array. Phy. Rev. Appl. (2019) doi: 10.1103/PhysRevApplied.13.031001

28. Huang, J., Zhou, X.: A time-varying mass metamaterial for non-reciprocal wave propagation. Int. J. Solids Struct. (2019) doi: 10.1016/j.ijsolstr.2018.12.029

29. Attarzadeh, M., Callanan, J., Nouh, M.: Experimental observation of non-reciprocal waves in a resonant metamaterial beam. Phy. Rev. Appl. (2020) doi: 10.1103/PhysRevApplied.13.021001

30. Berry, M. V.: Quantal phase factors accompanying adiabatic changes. Proc. R. Soc. A (1984) doi: 10.1098/rspa.1984.0023

31. Chaunsali, R., Li, F., Yang, J.: Stress wave isolation by purely mechanical topological phononic crystals. Scientific Reports (2016) doi: 10.1038/srep30662

32. Chen, H., Yao, L. Y., Nassar, H., Huang, G. L.: Mechanical quantum Hall effect in time-modulated elastic materials. Phys. Rev. Appl. (2019) doi: 10.1103/PhysRevApplied.11.044029

33. Cullen, A.: A travelling-wave parametric amplifier. Nature, (1958) doi: 10.1038/181332a0

34. Hayrapetyan, A., Grigoryan, K., Petrosyan, R., Fritzsche, S.: Propagation of sound waves through a spatially homogeneous but smoothly time-dependent medium. Ann. Phys. (2013) doi: 10.1016/j.aop.2013.02.014

35. Lurie, K. A., Weekes, S. L.: Wave propagation and energy exchange in a spatio-temporal material composite with rectangular microstructure. J. Math. Anal. Appl. (2006) doi: 10.1016/j.jmaa.2005.03.093

36. Milton, G. W., Mattei, O.: Field patterns: a new mathematical object. Proc. R. Soc. A (2017) doi: 10.1098/rspa.2016.0819

37. Torrent, D., Parnell, W. J., Norris, A. N.: Loss compensation in time-dependent elastic meta-
materials. Phys. Rev. B (2018) doi: 10.1103/physrevb.97.014105